\documentclass[12pt,a4paper]{amsart}
\usepackage{amsmath,amssymb,amsthm,mathtools}
\usepackage[margin=1in]{geometry}
\usepackage{booktabs,microtype,xcolor,array,needspace}
\usepackage{tikz,float}
\usepackage[colorlinks=true,linkcolor=blue!45!black,citecolor=blue!45!black,urlcolor=blue!45!black]{hyperref}
\numberwithin{equation}{section}
\newcommand{\ev}[1]{\left.#1\right|_{u=1,\,v=y}}
\newcommand{\bd}{\operatorname{bd}}
\newcommand{\corn}{\operatorname{corn}}
\newcommand{\ind}{\mathbf{1}}
\newtheorem{theorem}{Theorem}[section]
\newtheorem{proposition}[theorem]{Proposition}
\newtheorem{lemma}[theorem]{Lemma}
\newtheorem{corollary}[theorem]{Corollary}
\theoremstyle{definition}
\newtheorem{example}[theorem]{Example}
\theoremstyle{plain}
\title{A Murnaghan--Nakayama Rule for Dual Stable Grothendieck Functions}
\author{Dang Tuan Hiep}
\address{Da Lat University, 1 Phu Dong Thien Vuong, Lam Vien-Da Lat, Lam Dong, Viet Nam}
\email{hiepdt@dlu.edu.vn}
\author{Do Le Hai Thuy}
\address{Institute of Mathematics, Vietnam Academy of Science and Technology, 18 Hoang Quoc Viet, Hanoi, Viet Nam}
\email{cbl.dolehaithuy@gmail.com}
\date{November 3, 2026}
\begin{document}
\begin{abstract}
We give a closed boundary formula for multiplication by power sums in the dual stable Grothendieck basis. The off-diagonal coefficients are signed binomial coefficients supported on corner-connected broken ribbons, which are skew diagrams connected through edges or vertices and containing no square of four boxes. For such a diagram, the smallest power-sum index with a nonzero coefficient is the number of boxes plus the number of edge-connected components minus one. All larger indices give nonzero coefficients when the boundary statistic is positive and zero coefficients when it vanishes. The diagonal coefficient is minus the number of removable corners of the indexing partition. The formula evaluates the stable specialization of Nguyen's Murnaghan--Nakayama rule for K-k-Schur functions. We give two independent proofs, one using Newton's identity and the dual Pieri rules and the other using determinantal formulas and Schur complements.
\end{abstract}
\maketitle

\section{Introduction}
The classical Murnaghan--Nakayama rule expresses multiplication by a power sum in the Schur basis as a signed sum over ribbons, with signs determined by their heights \cite{Macdonald}. Stable Grothendieck functions and their Hall-dual basis are inhomogeneous $K$-theoretic analogues of Schur functions \cite{Buch,LP}. Lam and Pylyavskyy \cite{LP} describe the dual stable Grothendieck functions $g_\lambda$ as generating functions of reverse plane partitions, where the exponent of $x_i$ records the number of columns containing $i$. In this paper, we determine the coefficients
\[
 F_k(\lambda,\mu)=[g_\mu](g_\lambda p_k),\qquad k\ge1,
\]
by closed formulas depending only on the boundary geometry of $\lambda$ and $\mu$. The notation $[g_\mu](g_\lambda p_k)$ denotes the coefficient of $g_\mu$ in the expansion of $g_\lambda p_k$ in the dual stable Grothendieck basis.

Nguyen--Hiep--Son--Thuy \cite{NguyenEtAl} establish a Murnaghan--Nakayama rule for stable Grothendieck functions $G_\lambda$ of Grassmannian type. Kundu \cite{Kundu} gives an alternative proof based on alternants. The extension to the canonical stable Grothendieck functions $G_\lambda^{(\alpha,\beta)}$ uses the deformed power sum \cite[Section~3.3, Definition~3 and Theorem~2]{KunduCanonical}
\[
 p_k^\alpha(X)=p_k\left(\frac{x_1}{1-\alpha x_1},\frac{x_2}{1-\alpha x_2},\ldots\right).
\]
Hall duality relates power-sum multiplication in the $g$-basis to the action of the adjoint operator on the $G$-basis:
\[
 F_k(\lambda,\mu)=\langle G_\mu,g_\lambda p_k\rangle
 =\langle p_k^\perp G_\mu,g_\lambda\rangle,
\]
where $p_k^\perp$ denotes the Hall-adjoint of multiplication by $p_k$.

In the affine setting, Lam--Schilling--Shimozono \cite{LSS} realize the $K$-$k$-Schur functions as Schubert representatives for the $K$-homology of the affine Grassmannian of type~$A$. Bandlow--Schilling--Zabrocki \cite{BSZ} prove a Murnaghan--Nakayama rule for $k$-Schur functions using the affine nilCoxeter algebra. Nguyen \cite{NguyenKk} extends the rule to $K$-$k$-Schur functions, describing the coefficients by signed, weighted counts of weak hook words and providing both a formulation in terms of partitions and algorithms for their computation.

To compare this rule with ours, we reserve $k$ for the power-sum index and denote the affine parameter by $K$. Write $g_\rho^{(K)}$ for the corresponding member of the $K$-$k$-Schur family. A partition $\lambda$ is $K$-bounded if $\lambda_1\le K$. For such a partition, Nguyen's rule \cite[Corollaries~9.2 and~9.6]{NguyenKk} takes the form
\begin{equation}\label{eq:affine-expansion}
 p_k g_\lambda^{(K)}
 =\sum_{\substack{\mu_1\le K\\|\mu|\le |\lambda|+k}}
 C_k^{(K)}(\lambda,\mu)g_\mu^{(K)},\qquad 1\le k\le K.
\end{equation}
Condition (K2) of Nguyen's partition formulation \cite[Corollary~9.6]{NguyenKk} applies to the skew diagram between the $(K+1)$-cores associated with $\mu$ and $\lambda$. Nguyen calls a skew diagram a ribbon when it contains no $2\times2$ block \cite[Definitions~3.1]{NguyenKk}. For nonempty diagrams, this is our broken ribbon convention.

Blasiak--Morse--Seelinger \cite[Corollary~2.9]{BMS} prove the stability identity
\begin{equation}\label{eq:BMS-stability}
 g_\rho^{(K)}=g_\rho
 \quad\text{whenever}\quad \rho_1+\ell(\rho)-1\le K.
\end{equation}
Set $N=|\lambda|+k$ and choose $K\ge N$. Every nonempty partition $\rho$ with $|\rho|\le N$ satisfies
\[
 \rho_1+\ell(\rho)-1\le |\rho|\le N\le K,
\]
while $g_\varnothing^{(K)}=g_\varnothing=1$. Every hook length of such a partition $\rho$ is at most $K$, so the core-to-bounded-partition map retains every box and $\rho$ is its own associated $(K+1)$-core. Hence \eqref{eq:BMS-stability} applies to $\lambda$ and to every index occurring in \eqref{eq:affine-expansion}. Substituting these identities and comparing coefficients in the $g$-basis yields
\begin{equation}\label{eq:stable-identification}
 C_k^{(K)}(\lambda,\mu)=F_k(\lambda,\mu)
 \qquad\bigl(K\ge |\lambda|+k,\ |\mu|\le |\lambda|+k\bigr).
\end{equation}
Thus an explicit power-sum multiplication rule for dual stable Grothendieck functions already follows from Nguyen's affine rule and the stability theorem.

Our main contribution is a closed boundary formula for the stabilized coefficients. For a nonempty skew diagram $\theta=\mu/\lambda$, write $d=|\theta|$, let $r$ and $c$ denote its numbers of occupied rows and edge-connected components, respectively, and set $B=\bd(\lambda,\mu)$ as in \eqref{eq:boundary}. A corner-connected broken ribbon is a diagram without a $2\times2$ block that is connected under edge-or-vertex adjacency. For such a diagram in the range of \eqref{eq:stable-identification}, Theorem~\ref{thm:main} gives
\begin{equation}\label{eq:stable-boundary}
 C_k^{(K)}(\lambda,\mu)=F_k(\lambda,\mu)
 =(-1)^{r-c}\binom{k-d+B-c}{k-d-c+1}.
\end{equation}
All other off-diagonal coefficients vanish, while the diagonal coefficient is the negative number of removable corners of $\lambda$. For an allowed diagram, Corollary~\ref{cor:exact} gives the exact nonvanishing range: if $B>0$, the coefficient is nonzero precisely when $k\ge d+c-1$; if $B=0$, it is nonzero only when $k=d$. The formula therefore determines the support, minimal degree defect, and coefficient growth directly from boundary data, without enumerating weak hook words. Identity~\eqref{eq:stable-identification} identifies the aggregate coefficients; a direct correspondence between individual weak hook words and surviving boundary configurations is a separate combinatorial question.

Several related developments provide context for our methods. Galashin--Grinberg--Liu \cite{GGL} introduce refined dual stable Grothendieck polynomials and prove their symmetry using generalized Bender--Knuth involutions. Kim \cite{Kim} proves Grinberg's conjectured Jacobi--Trudi formula for these polynomials. Iwao \cite{Iwao} uses free fermions to derive determinantal formulas for skew stable Grothendieck polynomials and their duals. In the proof of Proposition~8.1 of the arXiv version, Iwao expresses multiplication of $g_{\lambda/\eta}$ by $p_k$ through a noncommutative supersymmetric power-sum operator acting on the pair $(\lambda,\eta)$. Taking $\eta=\varnothing$ gives an operator expression for $p_k g_\lambda$. The boundary formula above evaluates the coefficients of this product directly in the $g$-basis. Yeliussizov's work on deformations and skew Cauchy identities \cite{Yeliussizov,Yskew} provides the duality and Pieri formulas used in our proofs.

We present two independent proofs of the boundary formula. The first combines Newton's identity with the dual Pieri rules \cite{Yeliussizov}. The intermediate partitions are parametrized by independent row choices, so the associated generating polynomial factors over the rows. A first-order differential operator annihilates the contribution of every diagram with two or more components under edge-or-vertex adjacency. For a corner-connected broken ribbon, the sole surviving contribution comes from differentiating the unique vanishing terminal-row factor, yielding the boundary formula.

The second proof starts from the skew Jacobi--Trudi and bialternant formulas of Amanov--Yeliussizov \cite[Theorems~16 and~19]{AY}. Cauchy--Binet yields mutually inverse transition matrices between the $g$-basis and the Schur basis. Multiplication by a power sum is expressed as a sum of alternants in which one exponent is shifted, leading to a cofactor formula. Schur complements eliminate the unchanged rows and columns, and a rank drop forces vanishing for the forbidden diagrams. For an allowed diagram, the unique surviving cofactor yields the boundary formula. We also compute the diagonal correction $n-1$ arising when the finite convolution is extended to the full kernel. The determinant identities and exponent-shift principle are established tools; this proof uses them to obtain the closed boundary evaluation.

The two proofs give independent, complete derivations of the coefficient formula. Section~\ref{sec:notation} introduces the notation and states the main result. Sections~\ref{sec:rowproof} and~\ref{sec:detproof} present the Newton--Pieri and determinantal proofs, respectively.

\section{Notation and the main result}\label{sec:notation}
For a partition $\lambda$, write $|\lambda|$ for its size and $\ell(\lambda)$ for its length, and set $\lambda_i=0$ for $i>\ell(\lambda)$. We draw Young diagrams in English notation and write $\lambda\subseteq\mu$ when $\lambda_i\le\mu_i$ for every $i$. The number of removable corners of $\lambda$ is
\[
 \corn(\lambda)=\#\{i\ge1:\lambda_i>\lambda_{i+1}\}.
\]
Let
\[\Lambda=\bigoplus_{m\ge0}\Lambda_m\] 
be the graded ring of symmetric functions over $\mathbb Z$, where $\Lambda_m$ is its homogeneous component of degree $m$, and let
\[\widehat\Lambda=\prod_{m\ge0}\Lambda_m\] 
be its degree completion. Denote the complete homogeneous, elementary, power-sum, and Schur functions by $h_m,e_m,p_m,s_\lambda$, respectively. We adopt the conventions $h_0=e_0=1$ and $h_m=e_m=0$ for $m<0$.

The Hall inner product makes the Schur basis orthonormal and extends to a pairing $\widehat\Lambda\times\Lambda\to\mathbb Z$. The stable Grothendieck functions $G_\lambda\in\widehat\Lambda$ and their dual functions $g_\mu\in\Lambda$ satisfy
\[
 \langle G_\lambda,g_\mu\rangle=\delta_{\lambda\mu}.
\]
A reverse plane partition of a skew shape $\alpha/\beta$ is a filling by positive integers weakly increasing along rows from left to right and along columns from top to bottom. Write $\operatorname{RPP}(\alpha/\beta)$ for the set of such fillings and $\operatorname{col}_j(T)$ for the number of columns of $T$ containing $j$. The skew dual stable Grothendieck function is
\[
 g_{\alpha/\beta}(X)
 =\sum_{T\in\operatorname{RPP}(\alpha/\beta)}
   \prod_{j\ge1}x_j^{\operatorname{col}_j(T)},
 \qquad X=(x_1,x_2,\ldots),
\]
where $\beta\subseteq\alpha$ \cite{LP,AY}. In particular, $g_\alpha=g_{\alpha/\varnothing}$, $g_{\alpha/\alpha}=1$, and $g_{\alpha/\beta}(0)=0$ for $\beta\subsetneq\alpha$. The top homogeneous component of $g_\lambda$ is $s_\lambda$.

For $f\in\Lambda$, let $[g_\mu]f$ denote the coefficient of $g_\mu$ in its expansion in the basis $\{g_\rho\}$. We study
\[
 F_k(\lambda,\mu):=[g_\mu](g_\lambda p_k),\qquad k\ge1.
\]
A broken ribbon is a nonempty skew diagram containing no $2\times2$ block. It is a ribbon if it is connected by adjacency along a full edge, and a corner-connected broken ribbon if it is connected by adjacency along an edge or a vertex. The height of a ribbon $\theta$ is $r(\theta)-1$, where $r(\theta)$ is its number of occupied rows. For $\lambda\subsetneq\mu$, set $\theta=\mu/\lambda$ and write
\[
 d=|\theta|,\qquad r=r(\theta)=\#\{i:\mu_i>\lambda_i\},
 \qquad c=\kappa(\theta),
\]
where $\kappa(\theta)$ denotes the number of components under edge adjacency. Thus a corner-connected broken ribbon may have $c>1$. Define the boundary statistic
\begin{equation}\label{eq:boundary}
 B=\bd(\lambda,\mu)
 =\sum_{i\ge1}\bigl(\ind_{\mu_i>\lambda_{i+1}}-
                         \ind_{\lambda_i>\mu_{i+1}}\bigr).
\end{equation}
Here $\ind_P$ equals $1$ if $P$ holds and $0$ otherwise, and the sum is finite. The statistic $B$ depends on the indexing pair $(\lambda,\mu)$, including unchanged rows adjacent to the skew diagram. For example, $(3,2)/(3,1)$ and $(2,2)/(2,1)$ both consist of the single cell $(2,2)$, but their values of $B$ are $0$ and $1$, respectively. For a coefficient indexed by $\lambda\subseteq\mu$, its degree defect is $|\lambda|+k-|\mu|$; for $\lambda\subsetneq\mu$, this equals $k-d$.

We use generalized binomial coefficients, defined by the usual falling-factorial formula for nonnegative lower indices and by $\binom{x}{m}=0$ for $m<0$. In particular, $\binom{-1}{0}=1$. All generating series are interpreted as formal expansions at zero. For a formal Laurent series $f(z)$, the notation $[z^j]f(z)$ denotes the coefficient of $z^j$. If $f(z)$ is a power series, this coefficient is zero for $j<0$.

\Needspace{18\baselineskip}
\begin{samepage}
\begin{theorem}\label{thm:main}
Let $\lambda$ and $\mu$ be partitions. If $\lambda\nsubseteq\mu$, then $F_k(\lambda,\mu)=0$ for every $k\ge1$. On the diagonal,
\begin{equation}\label{eq:diagonal}
 F_k(\lambda,\lambda)=-\corn(\lambda),\qquad k\ge1.
\end{equation}
For $\lambda\subsetneq\mu$, the coefficients vanish for all $k\ge1$ unless $\mu/\lambda$ is a corner-connected broken ribbon. For such a diagram, $B\ge0$ and
\begin{align}
 \sum_{k\ge1}F_k(\lambda,\mu)z^k
 &=(-1)^{r-c}\frac{z^{d+c-1}}{(1-z)^B},\label{eq:gf}\\
 F_k(\lambda,\mu)
 &=(-1)^{r-c}\binom{k-d+B-c}{k-d-c+1},\qquad k\ge1.\label{eq:closed}
\end{align}
\end{theorem}
\end{samepage}

For $\lambda\subsetneq\mu$, put $d=|\mu/\lambda|$ and define the shifted series
\[
 \Phi_{\lambda,\mu}(t):=\sum_{m\ge0}F_{d+m}(\lambda,\mu)t^m.
\]
For a corner-connected broken ribbon, this series satisfies
\begin{equation}\label{eq:shifted}
 \Phi_{\lambda,\mu}(t)=(-1)^{r-c}t^{c-1}(1-t)^{-B},
\end{equation}
and it vanishes otherwise. Since $F_k(\lambda,\mu)=0$ for $k<d$, the unshifted series in \eqref{eq:gf} equals $z^d\Phi_{\lambda,\mu}(z)$. Thus the two proofs use the same generating series, with a shift of $d$ in the index.

\begin{lemma}\label{lem:support}
Let $f\in\Lambda$ and suppose that $[g_\mu](g_\lambda f)\ne0$. Then $\lambda\subseteq\mu$. If $\deg f\le k$ for an integer $k\ge0$, then also $|\mu|\le|\lambda|+k$.
\end{lemma}
\begin{proof}
Let $\Delta_\Lambda$ denote the coproduct, defined for disjoint alphabets by
\[\Delta_\Lambda f(X,Y)=f(X\cup Y)\] 
and extended continuously to the degree completion. The branching identity \cite[Proposition~4.3]{Yskew} gives
\[
 \Delta_\Lambda G_\mu=\sum_{\nu\subseteq\mu}G_\nu\otimes G_{\mu//\nu}.
\]
Here $G_{\mu//\nu}\in\widehat\Lambda$ is the skew stable Grothendieck function on $\mu/\nu$ together with the removable corners of $\nu$. In its set-valued tableau model, cells of $\mu/\nu$ must carry nonempty sets, whereas those removable corners may carry empty sets; this boundary condition accounts for the double slash. Compatibility of the Hall pairing with the coproduct yields
\begin{align*}
 [g_\mu](g_\lambda f)
 &=\langle G_\mu,g_\lambda f\rangle\\
 &=\langle\Delta_\Lambda G_\mu,g_\lambda\otimes f\rangle\\
 &=\begin{cases}
 \langle G_{\mu//\lambda},f\rangle,&\lambda\subseteq\mu,\\
 0,&\text{otherwise}.
 \end{cases}
\end{align*}
These pairings are well defined because $g_\lambda$ and $f$ have bounded degree. This proves containment.

For the degree bound, the unitriangular relations $g_\rho=s_\rho+$ terms of lower degree imply that a symmetric function of degree at most $N$ expands only in functions $g_\rho$ with $|\rho|\le N$. Since $\deg(g_\lambda f)\le|\lambda|+k$, the assertion follows.
\end{proof}

\Needspace{19\baselineskip}
\begin{samepage}
\begin{corollary}\label{cor:exact}
Let $\lambda\subsetneq\mu$. The coefficient $F_k(\lambda,\mu)$ vanishes unless $\mu/\lambda$ is a corner-connected broken ribbon. On this class, if $B>0$, then
\[
 F_k(\lambda,\mu)\ne0
 \quad\Longleftrightarrow\quad k\ge d+c-1.
\]
If $B=0$, then necessarily $r=c=1$, and
\[
 F_k(\lambda,\mu)\ne0\quad\Longleftrightarrow\quad k=d.
\]
\end{corollary}
\begin{proof}
The support restriction follows from Theorem~\ref{thm:main}. If $B>0$, every coefficient of $(1-z)^{-B}$ in nonnegative degree is strictly positive, so \eqref{eq:gf} gives the stated threshold. If $B=0$, the generating function is the single monomial $(-1)^{r-c}z^{d+c-1}$. For a corner-connected broken ribbon, the $r-1$ junctions between consecutive occupied rows each contribute $1$ to \eqref{eq:boundary}, and all other summands are nonnegative. Hence $B\ge r-1$, so $B=0$ forces $r=1$ and therefore $c=1$. Its only nonzero coefficient therefore occurs at $k=d$.
\end{proof}
\end{samepage}

\begin{corollary}
When $|\mu|-|\lambda|=k\ge1$, the coefficient $F_k(\lambda,\mu)$ is nonzero precisely when $\lambda\subseteq\mu$ and $\mu/\lambda$ is a ribbon. Its value is then $(-1)^{r-1}$.

For a fixed pair $(\lambda,\mu)$ such that $\mu/\lambda$ is a corner-connected broken ribbon and $B>0$, the coefficients $F_k(\lambda,\mu)$ for $k\ge d+c-1$ are given by a polynomial in $k$ of degree $B-1$, with leading coefficient
\[
 \frac{(-1)^{r-c}}{(B-1)!}.
\]
\end{corollary}
\begin{proof}
When $d=k$, \eqref{eq:gf} gives zero if $c>1$ and $(-1)^{r-1}$ if $c=1$. Together with the support restriction in Theorem~\ref{thm:main}, this proves the first assertion.

For $B>0$ and $k\ge d+c-1$, rewrite \eqref{eq:closed} as
\[
 F_k(\lambda,\mu)=(-1)^{r-c}\binom{k-d+B-c}{B-1}.
\]
The polynomial assertion and its leading coefficient follow.
\end{proof}

\begin{corollary}
For $k\ge1$,
\begin{equation}\label{eq:hooks}
 p_k=\sum_{\substack{a\ge1,\ b\ge0\\a+b\le k}}
 (-1)^b\binom{k-a}{b}\,g_{(a,1^b)}.
\end{equation}
\end{corollary}
\begin{proof}
For $\lambda=\varnothing$, the allowed nonempty straight shapes are exactly the hooks. If $\mu=(a,1^b)$, then
\[
 d=a+b,\qquad r=b+1,\qquad c=1,\qquad B=b+1.
\]
The degree bound restricts the sum to $a+b\le k$, and substitution into \eqref{eq:closed} gives the coefficient in \eqref{eq:hooks}. The empty partition contributes zero because $\corn(\varnothing)=0$.
\end{proof}

\begin{example}\label{ex:corner-ribbon}
Let $\lambda=(3,1,1)$ and $\mu=(4,3,2)$. The skew diagram $\theta=\mu/\lambda$, shown in Figure~\ref{fig:corner-ribbon}, has no $2\times2$ block. Its two edge-connected components meet at a single vertex, so it is a corner-connected broken ribbon with
\[
 d=4,\qquad r=3,\qquad c=2.
\]
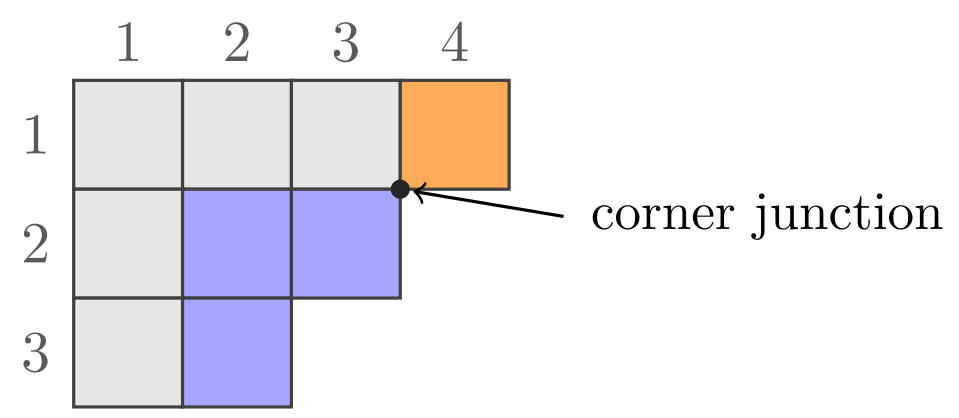
\begin{figure}[H]
\centering
\begin{tikzpicture}[x=0.72cm,y=0.72cm,font=\small,
 cell/.style={draw=black!75,line width=0.6pt},
 old/.style={cell,fill=gray!20},
 first/.style={cell,fill=orange!65},
 second/.style={cell,fill=blue!35}]
\foreach \x/\y in {0/0,1/0,2/0,0/-1,0/-2}
  \draw[old] (\x,\y) rectangle ++(1,-1);
\draw[first] (3,0) rectangle ++(1,-1);
\foreach \x/\y in {1/-1,2/-1,1/-2}
  \draw[second] (\x,\y) rectangle ++(1,-1);
\node at (3.5,-0.5) {};
\node at (1.5,-1.5) {};
\node at (2.5,-1.5) {};
\node at (1.5,-2.5) {};
\foreach \x/\j in {0.5/1,1.5/2,2.5/3,3.5/4}
  \node[text=black!65] at (\x,0.35) {\j};
\foreach \y/\i in {-0.5/1,-1.5/2,-2.5/3}
  \node[text=black!65] at (-0.35,\y) {\i};
\fill[black!85] (3,-1) circle (1.8pt);
\draw[->,line width=0.6pt] (4.5,-1.25) -- (3.12,-1.02);
\node[anchor=west,align=left,font=\footnotesize] at (4.55,-1.25)
  {corner\ junction};
\end{tikzpicture}
\caption{The diagram $(4,3,2)/(3,1,1)$. Gray cells form $\lambda$; the two colored edge-connected components meet at the marked vertex.}
\label{fig:corner-ribbon}
\end{figure}
Using \eqref{eq:boundary}, with zero parts appended, gives
\[
\begin{aligned}
 B&=(\ind_{4>1}-\ind_{3>3})
   +(\ind_{3>1}-\ind_{1>2})
   +(\ind_{2>0}-\ind_{1>0})\\
  &=1+1+0=2.
\end{aligned}
\]
Theorem~\ref{thm:main} therefore yields
\[
 \sum_{k\ge1}F_k\bigl((3,1,1),(4,3,2)\bigr)z^k
 =-\frac{z^5}{(1-z)^2}
 =-z^5-2z^6-3z^7-\cdots,
\]
or, equivalently,
\[
 F_k\bigl((3,1,1),(4,3,2)\bigr)
 =\begin{cases}
 0,&1\le k\le4,\\
 -(k-4),&k\ge5.
 \end{cases}
\]
Thus the first nonzero coefficient occurs at $k=d+c-1=5$, with minimal degree defect $c-1=1$. In particular, the coefficient at $k=d=4$ vanishes because $\theta$ is not edge-connected. The same partition $\lambda$ has two removable corners, so the diagonal coefficient is
\[
 [g_{(3,1,1)}]\bigl(g_{(3,1,1)}p_k\bigr)=-2
 \qquad(k\ge1).
\]
\end{example}

The row proof derives these support restrictions directly from the Pieri rules. The determinantal proof obtains them independently from Lemma~\ref{lem:support}.

\begin{example}\label{ex:two-one-p-three}
Let $\lambda=(2,1)$ and $k=3$. Lemma~\ref{lem:support} restricts the expansion to $\mu\supseteq\lambda$ with $|\mu|\le 6$. The two removable corners of $\lambda$ are $(1,2)$ and $(2,1)$, so the diagonal coefficient is $-2$ by \eqref{eq:diagonal}.

For $\mu\supsetneq\lambda$, put $d=|\mu/\lambda|$ and $j=3-d-c+1$, with $r$, $c$, and $B$ as in Theorem~\ref{thm:main}. Since $d\le3$, a $2\times2$ block cannot occur. For a corner-connected diagram, \eqref{eq:closed} gives
\[
 F_3((2,1),\mu)=(-1)^{r-c}\binom{B+j-1}{j}.
\]
The coefficient vanishes if $j<0$ or if the diagram is disconnected under edge-or-vertex adjacency.

For $\mu=(2,1,1)$, the only added box is $(3,1)$, so $d=r=c=1$. The three potentially nonzero summands of \eqref{eq:boundary} give
\[
 B=(\ind_{2>1}-\ind_{2>1})+(\ind_{1>0}-\ind_{1>1})
   +(\ind_{1>0}-\ind_{0>0})=2.
\]
Thus $j=2$ and the coefficient is $\binom{3}{2}=3$. For $\mu=(3,1)$, a single box is added in the first row. Here $B=0$ and the coefficient is $\binom{1}{2}=0$.

For $\mu=(3,2)$, the added boxes $(1,3)$ and $(2,2)$ share a vertex, so $d=r=c=2$, $B=1$, and $j=0$. The coefficient is therefore $1$. For $\mu=(3,2,1)$, the additional box $(3,1)$ extends this chain to three edge-connected components. The smallest index with a nonzero coefficient is $d+c-1=5$, so the coefficient at $k=3$ vanishes.

The table lists all off-diagonal candidates. A dash indicates that the diagram fails the edge-or-vertex connectivity condition.
\begin{center}
\small
\begin{tabular}{@{}r l r r r r r@{}}
\toprule
$d$ & $\mu$ & $r$ & $c$ & $B$ & $j$ & $F_3((2,1),\mu)$\\
\midrule
1 & $(3,1)$         & 1 & 1 & 0 &  2 & 0\\
1 & $(2,2)$         & 1 & 1 & 1 &  2 & 1\\
1 & $(2,1,1)$       & 1 & 1 & 2 &  2 & 3\\
\addlinespace
2 & $(4,1)$         & 1 & 1 & 0 &  1 & 0\\
2 & $(3,2)$         & 2 & 2 & 1 &  0 & 1\\
2 & $(3,1,1)$       & 2 & --- & --- & --- & 0\\
2 & $(2,2,1)$       & 2 & 2 & 3 &  0 & 1\\
2 & $(2,1,1,1)$     & 2 & 1 & 3 &  1 & $-3$\\
\addlinespace
3 & $(5,1)$         & 1 & 1 & 0 &  0 & 1\\
3 & $(4,2)$         & 2 & 2 & 1 & $-1$ & 0\\
3 & $(4,1,1)$       & 2 & --- & --- & --- & 0\\
3 & $(3,3)$         & 2 & 1 & 1 &  0 & $-1$\\
3 & $(3,2,1)$       & 3 & 3 & 3 & $-2$ & 0\\
3 & $(3,1,1,1)$     & 3 & --- & --- & --- & 0\\
3 & $(2,2,2)$       & 2 & 1 & 3 &  0 & $-1$\\
3 & $(2,2,1,1)$     & 3 & 2 & 4 & $-1$ & 0\\
3 & $(2,1,1,1,1)$   & 3 & 1 & 4 &  0 & 1\\
\bottomrule
\end{tabular}
\end{center}
Grouping terms by the number of added boxes gives
\begin{equation}\label{eq:two-one-p-three}
\begin{aligned}
g_{(2,1)}p_3={}&-2g_{(2,1)}
  +g_{(2,2)}+3g_{(2,1,1)}\\
 &+g_{(3,2)}+g_{(2,2,1)}-3g_{(2,1,1,1)}\\
 &+g_{(5,1)}-g_{(3,3)}-g_{(2,2,2)}+g_{(2,1,1,1,1)}.
\end{aligned}
\end{equation}
When $d=3$, the inequality $j\ge0$ forces $c=1$. The four surviving diagrams are ribbons, and their signs $(-1)^{r-1}$ agree with the classical rule.
\end{example}

\section{First proof}\label{sec:rowproof}
\subsection{The two dual Pieri rules}
A horizontal strip is a skew diagram with at most one box in each column, and a vertical strip is a skew diagram with at most one box in each row. The empty skew diagram is both a horizontal strip and a vertical strip. We write $\lambda'$ for the conjugate of $\lambda$.

In Yeliussizov's two-parameter notation \cite{Yeliussizov}, $g_\lambda=g_\lambda^{(0,1)}$. For $j\ge1$, we have $g_{(j)}=h_j$ and
\begin{equation}\label{row:eq:column}
 g_{(1^j)}=\sum_{i=1}^j\binom{j-1}{i-1}e_i,\qquad
 e_j=\sum_{i=1}^j(-1)^{j-i}\binom{j-1}{i-1}g_{(1^i)}.
\end{equation}
The second identity follows from the first by binomial inversion.

\begin{proposition}\label{row:prop:pieri}
For every partition $\lambda$ and integer $k\ge1$,
\begin{align}
 g_\lambda h_k
 &=\sum_{\substack{\mu\supseteq\lambda,\ \mu/\lambda\text{ horizontal}\\|\mu/\lambda|\le k}}
 (-1)^{k-|\mu/\lambda|}
 \binom{\#\{i:\lambda_i>\mu_{i+1}\}}{k-|\mu/\lambda|}g_\mu,\label{row:eq:hpieri}\\
 g_\lambda e_k
 &=\sum_{\substack{\nu\supseteq\lambda,\ \nu/\lambda\text{ vertical}\\|\nu/\lambda|\le k}}
 (-1)^{k-|\nu/\lambda|}
 \binom{\#\{j:\lambda'_j>\nu'_{j+1}\}+k-1}{k-|\nu/\lambda|}g_\nu.
 \label{row:eq:epieri}
\end{align}
\end{proposition}
\begin{proof}
We specialize the formulas in \cite[Proposition~8.12, Type~2]{Yeliussizov}. For a horizontal strip $\mu/\lambda$, put
\[
 d=|\mu/\lambda|,\qquad R=\#\{i:\lambda_i>\mu_{i+1}\}.
\]
The coefficient of $g_\mu^{(\alpha,\beta)}$ in $g_\lambda^{(\alpha,\beta)}h_k$ is
\[
 q_{\mu/\lambda}^{(k)}(\alpha,\beta)
 =\sum_{\ell=\max\{1,d\}}^k
 (-(\alpha+\beta))^{\ell-d}\binom{R}{\ell-d}
 (-\alpha)^{k-\ell}\binom{k-1}{\ell-1}.
\]
At $(\alpha,\beta)=(0,1)$, only the term $\ell=k$ contributes, yielding \eqref{row:eq:hpieri}.

For \eqref{row:eq:epieri}, the same proposition conjugates the partitions and interchanges the parameters. Setting
\[
 d=|\nu/\lambda|,\qquad R=\#\{j:\lambda'_j>\nu'_{j+1}\},
\]
we obtain
\begin{align*}
 q_{\nu'/\lambda'}^{(k)}(1,0)
 &=(-1)^{k-d}\sum_{\ell=\max\{1,d\}}^k
   \binom{R}{\ell-d}\binom{k-1}{\ell-1}\\
 &=(-1)^{k-d}\binom{R+k-1}{k-d},
\end{align*}
where the last equality follows from Vandermonde's identity. If $d>k$, the summation range is empty, so such shapes do not contribute.
\end{proof}

\subsection{The Newton--Pieri reduction}\label{row:sec:newton}
We use Newton's identity in the form
\begin{equation}\label{row:eq:newton}
 p_k=\sum_{i=0}^{k-1}(-1)^i(k-i)e_i h_{k-i},\qquad k\ge1.
\end{equation}
Indeed, let $H(z)=\sum_{j\ge0}h_jz^j$ and $E(z)=\sum_{j\ge0}e_jz^j$. Since $H(z)E(-z)=1$, we have
\[
 \sum_{k\ge1}p_kz^{k-1}=\frac{H'(z)}{H(z)}=H'(z)E(-z).
\]
Extracting the coefficient of $z^{k-1}$ gives \eqref{row:eq:newton}.

Throughout this proof, $d,r,c,B$ retain their meanings from Section~\ref{sec:notation}. Fix $\lambda\subsetneq\mu$ and define
\begin{equation}\label{row:eq:D}
 \mathcal D(\lambda,\mu)=\{\nu:\lambda\subseteq\nu\subseteq\mu,
 \ \nu/\lambda\text{ vertical},\ \mu/\nu\text{ horizontal}\}.
\end{equation}
For $\nu\in\mathcal D(\lambda,\mu)$, set
\begin{align*}
 a_\nu&=|\nu/\lambda|,&b_\nu&=|\mu/\nu|,\\
 E_\nu&=\#\{i:\lambda_i>\lambda_{i+1},\ \nu_i=\lambda_i\},&
 S_\nu&=\#\{i:\nu_i>\mu_{i+1}\}.
\end{align*}
Thus $a_\nu+b_\nu=d$. Here $E_\nu$ counts removable corners of $\lambda$ with no box of $\nu$ immediately to their right, and $S_\nu$ counts removable corners of $\nu$ with no box of $\mu$ immediately below. When $\nu$ is clear from the context, we abbreviate these quantities by $a,b,E,S$.

\begin{lemma}\label{row:lem:corners}
For every $\nu\in\mathcal D(\lambda,\mu)$,
\[
 \#\{j:\lambda'_j>\nu'_{j+1}\}=E_\nu.
\]
\end{lemma}
\begin{proof}
Suppose that $\lambda'_j>\nu'_{j+1}$ and set $i=\lambda'_j$. Since $\lambda\subseteq\nu$, we have
\[
 \lambda_{i+1}<j\le\lambda_i\le\nu_i\le j.
\]
Hence $\lambda_i=\nu_i=j$, so $(i,j)$ is a removable corner of $\lambda$ counted by $E_\nu$.

Conversely, if $\lambda_i>\lambda_{i+1}$ and $\nu_i=\lambda_i=j$, then $\lambda'_j=i>\nu'_{j+1}$. These assignments give a bijection between the columns counted on the left and the removable corners counted by $E_\nu$, proving the identity.
\end{proof}

\begin{proposition}\label{row:prop:sum}
Let $\lambda\subsetneq\mu$ and set $d=|\mu/\lambda|$. For every $m\ge0$,
\begin{equation}\label{row:eq:sum}
\begin{split}
 F_{d+m}(\lambda,\mu)=(-1)^m\sum_{\nu\in\mathcal D(\lambda,\mu)}(-1)^{a_\nu}
 \bigg[&b_\nu\binom{S_\nu-a_\nu-E_\nu}{m}\\
 &+S_\nu\binom{S_\nu-a_\nu-E_\nu-1}{m-1}\bigg].
\end{split}
\end{equation}
Equivalently,
\begin{equation}\label{row:eq:Phi}
 \Phi_{\lambda,\mu}(t)=\sum_{\nu\in\mathcal D(\lambda,\mu)}
 (-1)^{a_\nu}(1-t)^{S_\nu-a_\nu-E_\nu}
 \left(b_\nu-\frac{t}{1-t}S_\nu\right).
\end{equation}
\end{proposition}
\begin{proof}
Fix $\nu\in\mathcal D(\lambda,\mu)$ and abbreviate $a=a_\nu$, $b=b_\nu$, $E=E_\nu$, and $S=S_\nu$, so that $a+b=d$. Set $k=d+m$. In the $i$th term of \eqref{row:eq:newton}, apply the column Pieri rule from $\lambda$ to $\nu$, followed by the row Pieri rule from $\nu$ to $\mu$. Writing $i=a+q$, the total contribution of $\nu$ is
\[
 (-1)^{m+a}\sum_{q=0}^m(-1)^q(b+m-q)
 \binom{a+E+q-1}{q}\binom{S}{m-q}.
\]
The Pieri degree restrictions give $0\le q\le m$. If $b=0$, the endpoint $q=m$ corresponds to $i=k$ and may be included because its factor $k-i$ vanishes. At $i=0$, the same expression remains valid by $e_0=1$ and our generalized binomial convention.

The formal binomial identity
\[
 \sum_{q\ge0}(-1)^q\binom{a+E+q-1}{q}x^q=(1+x)^{-a-E},
\]
together with differentiation of $(1+x)^S$ to account for the factor $m-q$, shows that the inner sum equals
\[
 b[x^m](1+x)^{S-a-E}+S[x^{m-1}](1+x)^{S-a-E-1}.
\]
Extracting coefficients proves \eqref{row:eq:sum}. Summing over $m\ge0$ amounts to substituting $x=-t$ in the corresponding generating function, yielding \eqref{row:eq:Phi}.
\end{proof}

\subsection{The row proof and the generating function}\label{row:sec:factor}
Set $L=\ell(\mu)$ and extend $\lambda$ and $\mu$ by zero parts. For $1\le i\le L$, define
\[
 \delta_i=\mu_i-\lambda_i,\qquad
 \Delta_i=\mu_{i+1}-\lambda_i,\qquad
 \eta_i=\ind_{\lambda_i>\lambda_{i+1}}.
\]
We call row $i$ active if $\delta_i>0$, equivalently, if it contains a box of $\mu/\lambda$.

\begin{lemma}\label{row:lem:independent}
The assignments $\nu_i=\lambda_i+x_i$ give a bijection between $\mathcal D(\lambda,\mu)$ and $\prod_{i=1}^L X_i$, where
\begin{equation}\label{row:eq:Xi}
 X_i=\{x\in\{0,1\}:\mu_{i+1}\le\lambda_i+x\le\mu_i\}.
\end{equation}
In particular, $\mathcal D(\lambda,\mu)=\varnothing$ if and only if $\mu/\lambda$ contains a $2\times2$ block.
\end{lemma}
\begin{proof}
If $\nu\in\mathcal D(\lambda,\mu)$, the vertical-strip condition gives $\nu_i-\lambda_i\in\{0,1\}$, while the horizontal-strip condition is equivalent to the interlacing inequalities $\mu_i\ge\nu_i\ge\mu_{i+1}$. Hence $x_i=\nu_i-\lambda_i$ belongs to $X_i$.

Conversely, any choice $x_i\in X_i$ defines $\nu_i=\lambda_i+x_i$ satisfying
\[
 \nu_i\ge\mu_{i+1}\ge\nu_{i+1}.
\]
Thus $\nu$ is a partition belonging to $\mathcal D(\lambda,\mu)$, and the row choices are independent.

A $2\times2$ block across rows $i,i+1$ exists precisely when $\mu_{i+1}\ge\lambda_i+2$, in which case $X_i$ is empty. If no such block exists, then $\Delta_i\le1$ for every $i$. An inactive row has $X_i=\{0\}$. For an active row, $X_i=\{1\}$ when $\Delta_i=1$, and $X_i=\{0,1\}$ when $\Delta_i\le0$. Hence every $X_i$ is nonempty, proving the final assertion.
\end{proof}

Assume that $\theta=\mu/\lambda$ contains no $2\times2$ block. Let $u,v,y$ be independent indeterminates and define
\begin{align}
 \sigma_i(x)&=\ind_{\lambda_i+x>\mu_{i+1}},\qquad x\in\{0,1\},\label{row:eq:six}\\
 f_i(u,v;y)&=\sum_{x\in X_i}(-u)^xv^{\sigma_i(x)}y^{-x-\eta_i(1-x)},\label{row:eq:fi}\\
 Z(u,v;y)&=\prod_{i=1}^L f_i(u,v;y).\label{row:eq:Z}
\end{align}
Each $f_i$, and hence $Z$, lies in $\mathbb Z[y,y^{-1}][u,v]$. Substituting $y=1-t$ gives well-defined formal power series in $\mathbb Z[u,v][[t]]$.

Under the bijection of Lemma~\ref{row:lem:independent}, we have
\[
 a_\nu=\sum_i x_i,\qquad S_\nu=\sum_i\sigma_i(x_i),\qquad
 E_\nu=\sum_i\eta_i(1-x_i).
\]
Multiplying the row factors therefore gives
\[
 Z(u,v;y)=\sum_{\nu\in\mathcal D(\lambda,\mu)}
 (-1)^{a_\nu}u^{a_\nu}v^{S_\nu}y^{-a_\nu-E_\nu}.
\]
Using $b_\nu=d-a_\nu$, we can rewrite \eqref{row:eq:Phi} as
\begin{equation}\label{row:eq:differential}
 \Phi_{\lambda,\mu}(t)
 =\left.\left(d-u\frac{\partial}{\partial u}
 -(1-y)\frac{\partial}{\partial v}\right)Z(u,v;y)
 \right|_{u=1,\,v=y,\,y=1-t}.
\end{equation}
The partial derivatives are taken with $y$ held fixed, before the indicated substitutions.

\begin{lemma}\label{row:lem:factors}
For an inactive row, set
\[
 \tau_i=\eta_i-\sigma_i(0),\qquad \sigma_i(0)=\ind_{\lambda_i>\mu_{i+1}}.
\]
Then $\tau_i\in\{0,1\}$ and
\[
 f_i(u,v;y)=v^{\sigma_i(0)}y^{-\eta_i},\qquad f_i(1,y;y)=y^{-\tau_i}.
\]
For an active row, the factors and their specializations are
\begin{equation}\label{row:eq:table}
\begin{array}{c|c|c}
\text{Condition}&f_i(u,v;y)&f_i(1,y;y)\\\hline
\Delta_i=1&-u/y&-1/y\\
\Delta_i=0,\ \eta_i=1&(1-uv)/y&(1-y)/y\\
\Delta_i=0,\ \eta_i=0&1-uv/y&0\\
\Delta_i\le-1&v(1-u)/y&0
\end{array}
\end{equation}
\end{lemma}
\begin{proof}
For an inactive row, $X_i=\{0\}$. Since $\mu_{i+1}\ge\lambda_{i+1}$, we have $\sigma_i(0)\le\eta_i$. Thus $\tau_i\in\{0,1\}$, and the stated formulas follow directly from \eqref{row:eq:fi}.

Now suppose that row $i$ is active. If $\Delta_i=1$, then $X_i=\{1\}$ and $\sigma_i(1)=0$, giving $f_i=-u/y$. If $\Delta_i=0$, then $X_i=\{0,1\}$, $\sigma_i(0)=0$, and $\sigma_i(1)=1$, so
\[
 f_i=y^{-\eta_i}-uv/y.
\]
The two possible values of $\eta_i$ give the middle two rows of \eqref{row:eq:table}. Finally, if $\Delta_i\le-1$, then $X_i=\{0,1\}$ and $\sigma_i(0)=\sigma_i(1)=1$. Moreover, $\lambda_i>\mu_{i+1}\ge\lambda_{i+1}$ forces $\eta_i=1$, yielding $f_i=v(1-u)/y$. Specializing $u=1$ and $v=y$ gives the final column.
\end{proof}

\begin{lemma}\label{row:lem:zeros}
The factors $f_i(1,y;y)$ that vanish identically in $y$ are precisely those indexed by the last rows of the components of $\theta$ under edge-or-vertex adjacency. Hence their number equals the number of these components.
\end{lemma}
\begin{proof}
The boxes in each active row form a connected interval. Two adjacent active rows $i,i+1$ are joined under edge-or-vertex adjacency exactly when $\Delta_i\ge0$. Since $\theta$ contains no $2\times2$ block, $\Delta_i=1$ gives an edge junction and $\Delta_i=0$ gives a vertex junction.

If $\Delta_i=1$, the next row is necessarily active because $\mu_{i+1}=\lambda_i+1>\lambda_{i+1}$. If $\Delta_i=0$, the next row is active exactly when $\lambda_i>\lambda_{i+1}$, equivalently, when $\eta_i=1$. Thus, an active row is nonterminal precisely in the first two cases of \eqref{row:eq:table}, and terminal in the last two cases.

The table shows that the specialized factors are nonzero Laurent polynomials in the nonterminal cases and vanish identically in the terminal cases. Inactive rows also have nonzero factors by Lemma~\ref{row:lem:factors}. This proves the stated correspondence.
\end{proof}

\begin{proposition}\label{row:prop:gen}
The series $\Phi_{\lambda,\mu}(t)$ vanishes unless $\theta=\mu/\lambda$ is a corner-connected broken ribbon. In that case, \eqref{eq:shifted} holds.
\end{proposition}
\begin{proof}
If $\theta$ contains a $2\times2$ block, Lemma~\ref{row:lem:independent} gives $\mathcal D(\lambda,\mu)=\varnothing$, and \eqref{row:eq:Phi} implies $\Phi_{\lambda,\mu}=0$. Assume henceforth that no such block occurs. If $\theta$ has at least two components under edge-or-vertex adjacency, Lemma~\ref{row:lem:zeros} gives at least two factors of $Z$ vanishing at $u=1,v=y$. Since the operator in \eqref{row:eq:differential} differentiates at most one factor in each term, every resulting term retains a vanishing factor. Hence $\Phi_{\lambda,\mu}=0$.

Now suppose that $\theta$ is corner-connected. There is exactly one vanishing factor $f_j$, indexed by the last active row. Write
\[
 \mathcal R(u,v;y)=\prod_{\substack{1\le i\le L\\i\ne j}}f_i(u,v;y).
\]
The zeroth-order term and all terms obtained by differentiating $\mathcal R$ vanish after specialization. The product rule therefore gives
\begin{equation}\label{row:eq:terminal}
 \Phi_{\lambda,\mu}(t)
 =\left[-\ev{\left(u\frac{\partial f_j}{\partial u}
 +(1-y)\frac{\partial f_j}{\partial v}\right)}
 \mathcal R(1,y;y)\right]_{y=1-t}.
\end{equation}
Since $j$ is the last active row, $\mu_{j+1}=\lambda_{j+1}$. Set
\[
 \tau_+=\ind_{\lambda_j=\lambda_{j+1}}.
\]
If $\Delta_j\le-1$, then $f_j=v(1-u)/y$, and the differentiated contribution in \eqref{row:eq:terminal}, including its minus sign, equals $1$. If $\Delta_j=0$, then necessarily $\eta_j=0$, so $f_j=1-uv/y$ and this contribution equals $1/y$. In both cases, it is $y^{-\tau_+}$.

Let $n_E$ and $n_C$ count the edge and corner junctions, respectively, between consecutive active rows, and put
\[
 T=\sum_{\substack{1\le i\le L\\\delta_i=0}}\tau_i.
\]
Lemma~\ref{row:lem:factors} now yields
\begin{equation}\label{row:eq:preB}
 \Phi_{\lambda,\mu}(t)
 =(-1)^{n_E}t^{n_C}(1-t)^{-(n_E+n_C+T+\tau_+)}.
\end{equation}
Each active row is edge-connected, so $c=r-n_E$. Corner connectivity gives $n_E+n_C=r-1$, and consequently
\[
 n_E=r-c,\qquad n_C=c-1.
\]
It remains to identify the exponent of $1-t$. For an inactive row, the summand in \eqref{eq:boundary} is $\eta_i-\sigma_i(0)=\tau_i$. For an active row, its first indicator is $1$, so the summand equals $1$ if $\Delta_i\ge0$ and $0$ otherwise. Thus, every nonterminal active row contributes $1$, while the last row contributes $\tau_+$. Therefore
\begin{equation}\label{row:eq:Bidentity}
 B=n_E+n_C+T+\tau_+.
\end{equation}
Substituting these identities into \eqref{row:eq:preB} proves \eqref{eq:shifted}.
\end{proof}

\begin{proposition}\label{row:prop:diag}
For every $k\ge1$,
\[
 [g_\lambda](g_\lambda p_k)=-\corn(\lambda).
\]
\end{proposition}
\begin{proof}
For $\mu=\lambda$, the only intermediate partition is $\nu=\lambda$. Write $q=\corn(\lambda)$. Then $a_\nu=b_\nu=0$ and $E_\nu=S_\nu=q$. Repeating the Newton--Pieri coefficient extraction from the proof of Proposition~\ref{row:prop:sum}, with $d=0$ and $m=k\ge1$, gives
\[
 [g_\lambda](g_\lambda p_k)
 =(-1)^kq\binom{-1}{k-1}
 =(-1)^kq(-1)^{k-1}=-q.\qedhere
\]
\end{proof}

\begin{proof}[Completion of the proof of Theorem~\ref{thm:main}]
The Pieri rules and \eqref{row:eq:newton} imply that $F_k(\lambda,\mu)=0$ unless $\lambda\subseteq\mu$ and $|\mu/\lambda|\le k$. Proposition~\ref{row:prop:diag} gives the diagonal coefficient.

For $\lambda\subsetneq\mu$, Proposition~\ref{row:prop:gen} gives vanishing unless $\mu/\lambda$ is a corner-connected broken ribbon. On this class, \eqref{row:eq:Bidentity} shows that $B\ge0$. Substituting $t=z$ in \eqref{eq:shifted} and multiplying by $z^d$ proves \eqref{eq:gf}. Consequently,
\[
 F_k(\lambda,\mu)=(-1)^{r-c}[t^{k-d-c+1}](1-t)^{-B}.
\]
For $j\ge0$, the formal binomial theorem gives
\[
 [t^j](1-t)^{-B}=\binom{B+j-1}{j},
\]
including $B=0$ under our generalized binomial convention. For $j<0$, both sides vanish. Taking $j=k-d-c+1$ proves \eqref{eq:closed}.
\end{proof}

\begin{example}\label{ex:two-one-row-proof}
We recover several entries of Example~\ref{ex:two-one-p-three} directly from \eqref{row:eq:differential}. Fix $\lambda=(2,1)$ and write $\Phi_\mu=\Phi_{\lambda,\mu}$. The coefficient of $g_\mu$ is $[t^{3-d}]\Phi_\mu(t)$ off the diagonal.

For $\mu=(3,2)$, the choices are $X_1=X_2=\{0,1\}$. The four intermediate partitions are $(2,1)$, $(3,1)$, $(2,2)$, and $(3,2)$. The factors are
\[
 f_1=\frac{1-uv}{y},\qquad f_2=\frac{v(1-u)}{y}.
\]
At $u=1$ and $v=y$, only $f_2$ vanishes. Its differentiated contribution in \eqref{row:eq:terminal} is $1$, so $\Phi_{(3,2)}(t)=t/(1-t)$ and $[t]\Phi_{(3,2)}=1$.

For $\mu=(2,1,1)$, the factors are $v/y$, $1/y$, and $1-uv/y$. The last factor contributes $1/y$ after differentiation. Multiplying by the other two specialized factors gives $\Phi_{(2,1,1)}(t)=(1-t)^{-2}$, with $[t^2]\Phi_{(2,1,1)}=3$. For $\mu=(2,1,1,1)$, the first two factors are again $v/y$ and $1/y$, followed by $-u/y$ and $1-uv/y$. The last factor again contributes $1/y$, while the third specializes to $-1/y$. Hence $\Phi_{(2,1,1,1)}(t)=-(1-t)^{-3}$ and $[t]\Phi_{(2,1,1,1)}=-3$.

For $\mu=(3,1,1)$, the product is
\[
 Z(u,v;y)=\frac{v(1-u)}{y}\cdot\frac1y\cdot(1-uv/y).
\]
The first and third factors both vanish after specialization. Every term in \eqref{row:eq:differential} still has at least one of those vanishing factors, so $\Phi_{(3,1,1)}=0$. On the diagonal, $q=2$ in Proposition~\ref{row:prop:diag}, giving the coefficient $(-1)^3\,2\binom{-1}{2}=-2$.
\end{example}

\section{Second proof}\label{sec:detproof}
This proof is independent of Section~\ref{sec:rowproof} and uses the common preliminaries together with the formulas of Amanov--Yeliussizov. The boldface symbols $\mathbf u,\mathbf v$ denote vectors and are distinct from the scalar auxiliary variables $u,v$ used in the first proof. We denote the reduced matrix by $\mathcal L$; the symbol $L$ in the first proof denotes the number of rows.

\subsection{Determinants and inverse transition matrices}\label{sec:transition}
For integers $s,m$, define the scalar kernel
\begin{equation}\label{eq:H}
 \mathsf H_s(m)=[w^m](1-w)^{-s}
 =\begin{cases}
 \displaystyle\binom{s+m-1}{m},&m\ge0,\\
 0,&m<0.
 \end{cases}
\end{equation}
The symbol $\mathsf H$ distinguishes this kernel from the complete homogeneous symmetric functions $h_m$. In particular,
\[
 \mathsf H_0(m)=\delta_{m0},\qquad
 \mathsf H_{-q}(m)=(-1)^m\binom{q}{m}\quad(q\ge0).
\]
Let $\varphi$ be the automorphism of $\Lambda$ determined by
\[
 \varphi(h_m)=\sum_{j=0}^m h_j\qquad(m\ge1).
\]
For every integer $s$, its powers satisfy
\[
 \varphi^s(h_m)=\sum_{t\ge0}\mathsf H_s(t)h_{m-t}.
\]
For the alphabet $X=(x_1,x_2,\ldots)$, the formulas of Amanov--Yeliussizov \cite[Theorems~16 and~19]{AY} give
\begin{align}
 g_{\alpha/\beta}(X)
 &=\det\left[\varphi^{i-j}(h_{\alpha_i-\beta_j-i+j})(X)\right]_{i,j=1}^n,
 \label{eq:AYskew}\\
 g_\alpha(x_1,\ldots,x_n)
 &=\frac{\det\left[h_{\alpha_i+n-i}(1^{i-1},x_j)\right]_{i,j=1}^n}
 {\prod_{1\le i<j\le n}(x_i-x_j)}.\label{eq:AYalt}
\end{align}
Here $\beta\subseteq\alpha$ in \eqref{eq:AYskew}, and $n\ge\ell(\alpha)$ in both formulas, with partitions padded by zeros. The notation $1^{i-1}$ denotes an alphabet consisting of $i-1$ copies of $1$.

\begin{proposition}\label{prop:transition}
Let $\alpha,\rho$ be partitions, padded with zeros to a common length $n\ge\max\{\ell(\alpha),\ell(\rho)\}$, and set
\begin{align}
 A_{\alpha\rho}
 &=\det\left[\mathsf H_{i-1}(\alpha_i-\rho_j-i+j)\right]_{i,j=1}^n,\label{eq:A}\\
 D_{\alpha\rho}
 &=\det\left[\mathsf H_{1-j}(\alpha_i-\rho_j-i+j)\right]_{i,j=1}^n.\label{eq:D}
\end{align}
These coefficients are independent of $n$ once it includes both lengths, and
\begin{equation}\label{eq:change}
 g_\alpha=\sum_{\rho\subseteq\alpha}A_{\alpha\rho}s_\rho,
 \qquad s_\alpha=\sum_{\rho\subseteq\alpha}D_{\alpha\rho}g_\rho.
\end{equation}
\end{proposition}
\begin{proof}
Appending a zero part to both partitions adds a final row $(0,\ldots,0,1)$ to each matrix, so both determinants stabilize. Moreover, both vanish unless $\rho\subseteq\alpha$. Indeed, if $\rho_t>\alpha_t$, every entry with $i\ge t$ and $j\le t$ is zero. The first $t$ columns are therefore supported in the first $t-1$ rows and are linearly dependent.

Fix a size bound $M$ and take $n\ge\max\{1,M\}$. Specialization to $n$ variables is injective on symmetric functions of degree at most $M$. Expand each numerator entry in \eqref{eq:AYalt} as
\[
 h_{\alpha_i+n-i}(1^{i-1},x)
 =\sum_{t\ge0}\mathsf H_{i-1}(\alpha_i+n-i-t)x^t.
\]
Applying Cauchy--Binet with decreasing exponents $t_j=\rho_j+n-j$, and dividing by the Vandermonde, gives the first identity in \eqref{eq:change}.

To verify the inverse, put $a_i=\alpha_i-i$ and $b_j=\beta_j-j$. A second application of Cauchy--Binet gives
\begin{align*}
 \sum_\rho A_{\alpha\rho}D_{\rho\beta}
 &=\det\left[\sum_{t\ge-n}\mathsf H_{i-1}(a_i-t)
 \mathsf H_{1-j}(t-b_j)\right]_{i,j=1}^n\\
 &=\det\left[\mathsf H_{i-j}(a_i-b_j)\right]_{i,j=1}^n.
\end{align*}
For the second equality, the lower cutoff may be removed because the second factor forces $t\ge b_j\ge-n$. Scalar convolution then follows from $(1-w)^{-(i-1)}(1-w)^{-(1-j)}=(1-w)^{-(i-j)}$.

If $\beta\subseteq\alpha$, \eqref{eq:AYskew} identifies the last determinant with
\[
 g_{\alpha/\beta}(0)=\delta_{\alpha\beta},
\]
since a nonempty skew diagram has no constant term. If $\beta\nsubseteq\alpha$, the last determinant vanishes by the same zero-block argument, while every summand on the left vanishes by containment. Thus $AD=I$.

On the finite set of partitions of size at most $M$, both matrices are unitriangular, so $DA=I$ as well. Applying this inverse to the first identity proves the second. Since $M$ is arbitrary and the determinants stabilize, the result follows.
\end{proof}

\subsection{A cofactor formula with its cutoff correction}\label{sec:variation}
Fix $\lambda\subseteq\mu$ and $k\ge1$. Throughout the finite derivation, choose a single integer $n$, pad the partitions with zeros, and set
\begin{equation}\label{eq:nchoice}
 n\ge|\lambda|+k+\ell(\mu)+1,
 \qquad a_i=\lambda_i-i,\quad b_j=\mu_j-j.
\end{equation}
Define the $n\times n$ matrices $Q$ and $R^{(k)}$ by
\begin{equation}\label{eq:QR}
 Q_{ij}=\mathsf H_{i-j}(a_i-b_j),\qquad
 R^{(k)}_{ij}=\mathsf H_{i-j}(a_i-b_j+k).
\end{equation}
For any square matrices $Q,R$ of the same size, the linear term of the determinant satisfies
\begin{equation}\label{eq:variation}
 [\varepsilon]\det(Q+\varepsilon R)
 =\operatorname{tr}\bigl(\operatorname{adj}(Q)R\bigr)
 =\sum_{i,j}\operatorname{cof}_{ij}(Q)R_{ij}.
\end{equation}
Here
\[
 \operatorname{cof}_{ij}(Q)=(-1)^{i+j}\det Q^{(i\mid j)},\qquad
 (\operatorname{adj}Q)_{ji}=\operatorname{cof}_{ij}(Q),
\]
where $Q^{(i\mid j)}$ is obtained by deleting row $i$ and column $j$. The notation $\operatorname{tr}$ denotes the trace, and $[\varepsilon]$ extracts the coefficient of the formal parameter $\varepsilon$. Equivalently, \eqref{eq:variation} may be read in the ring of dual numbers, where $\varepsilon^2=0$. No invertibility assumption on $Q$ is needed.

\begin{proposition}\label{prop:cofactor}
With the notation and choice of $n$ in \eqref{eq:nchoice}--\eqref{eq:QR},
\begin{equation}\label{eq:cofactor}
 F_k(\lambda,\mu)=[\varepsilon]\det(Q+\varepsilon R^{(k)})
 -(n-1)\delta_{\lambda\mu}.
\end{equation}
\end{proposition}
\begin{proof}
We make the exponent-shift argument between the two transitions in \eqref{eq:change} explicit. Order the finite set
\[
 \mathcal E_n=\{a_1,a_1-1,\ldots,-n\}
\]
decreasingly, and define rectangular matrices, indexed by $1\le i,j\le n$ and $t\in\mathcal E_n$, by
\begin{align*}
 \mathsf U_{it}&=\mathsf H_{i-1}(a_i-t),\\
 \mathsf V_{tj}&=\mathsf H_{1-j}(t-b_j),\\
 \mathsf W^{(k)}_{tj}&=\mathsf H_{1-j}(t+k-b_j).
\end{align*}
The first kernel vanishes for $t>a_1$, while the second vanishes for $t<b_j$, where $b_j\ge-n$. Thus scalar convolution gives $\mathsf U\mathsf V=Q$ without any discarded terms.

Choose $n$ distinct exponents $t_1>\cdots>t_n$ from $\mathcal E_n$. They determine a partition $\rho$ by $t_i=\rho_i-i$, and the corresponding minor of $\mathsf U$ is $A_{\lambda\rho}$. The associated Schur alternant has nonnegative exponents $t_i+n=\rho_i+n-i$, which explains the lower cutoff $t_i\ge-n$. By Proposition~\ref{prop:transition}, this minor vanishes unless $\rho\subseteq\lambda$. Conversely, every such $\rho$ supplies a choice of exponents in $\mathcal E_n$. Cauchy--Binet therefore gives
\[
\begin{aligned}
 &\det\bigl(Q+\varepsilon\mathsf U\mathsf W^{(k)}\bigr)\\
 &\quad=\sum_{\rho\subseteq\lambda}A_{\lambda\rho}
 \det\left[\mathsf H_{1-j}(t_i-b_j)
       +\varepsilon\mathsf H_{1-j}(t_i+k-b_j)\right]_{i,j=1}^n,
 \qquad t_i=\rho_i-i.
\end{aligned}
\]
By multilinearity in the rows, the coefficient of $\varepsilon$ in each determinant on the right is the sum obtained by replacing exactly one exponent $t_a$ by $t_a+k$. If two resulting exponents coincide, the determinant has equal rows and is zero. Otherwise, let $\pi_a$ be the permutation sorting the shifted exponents as $\widehat t_1>\cdots>\widehat t_n$, and put $\sigma^{(a)}_i=\widehat t_i+i$. Sorting the rows shows that the contribution is
\[
 \operatorname{sgn}(\pi_a)
 \det\left[\mathsf H_{1-j}(\widehat t_i-b_j)\right]_{i,j=1}^n
 =\operatorname{sgn}(\pi_a)D_{\sigma^{(a)},\mu}.
\]
Here $\sigma^{(a)}$ is a partition of size $|\rho|+k$: the shifted exponents remain at least $-n$, and distinct decreasing integral exponents give weakly decreasing nonnegative parts after adding $i$.

On the other hand, multiplication of a Schur alternant by $p_k(x_1,\ldots,x_n)$ increases exactly one exponent by $k$. Hence
\[
 p_k s_\rho=\sum_{a=1}^n
 \operatorname{sgn}(\pi_a)s_{\sigma^{(a)}},
\]
where a term with repeated shifted exponents is interpreted as zero. Expanding $g_\lambda$ in the Schur basis and then each $s_{\sigma^{(a)}}$ in the $g$-basis using \eqref{eq:change}, we conclude that the coefficient of $\varepsilon$ in the Cauchy--Binet identity is $F_k(\lambda,\mu)$. Setting $\widetilde R^{(k)}=\mathsf U\mathsf W^{(k)}$ gives the exact finite formulas
\begin{align}
 F_k(\lambda,\mu)
 &=[\varepsilon]\det(Q+\varepsilon\widetilde R^{(k)}),\label{eq:finite}\\
 \widetilde R^{(k)}_{ij}
 &=\sum_{t\ge-n}\mathsf H_{i-1}(a_i-t)
 \mathsf H_{1-j}(t+k-b_j).\label{eq:Rfinite}
\end{align}
The sum in \eqref{eq:Rfinite} is finite because its first factor forces $t\le a_i\le a_1$. Every intermediate partition has size at most $|\lambda|+k$, and our choice of $n$ makes specialization to $n$ variables injective in this degree range. Thus the calculation recovers the stable coefficient. If $|\mu|>|\lambda|+k$, both sides of \eqref{eq:finite} vanish, by Lemma~\ref{lem:support} and containment for $D$, respectively. In particular, the lower cutoff $t\ge-n$ in \eqref{eq:Rfinite} is retained before any full convolution is used.

Extending the sum in \eqref{eq:Rfinite} to all integers $t$ gives $R^{(k)}$ by scalar convolution. Put
\[
 \mathcal T=R^{(k)}-\widetilde R^{(k)}.
\]
By \eqref{eq:variation}, it remains to prove
\[
 \sum_{i,j}\operatorname{cof}_{ij}(Q)\mathcal T_{ij}
 =(n-1)\delta_{\lambda\mu}.
\]
For a discarded term with $t<-n$, a nonzero second factor requires
\[
 t+k\ge b_j\ge-j.
\]
Hence $\mathcal T$ is supported in columns
\[
 j\ge n+1-k>\ell(\mu).
\]
The matrix $Q$ is upper triangular: if $i>j$, then $\lambda_i-\mu_j-i+j<0$. Its diagonal entry is $1$ when $\lambda_i=\mu_i$ and $0$ otherwise. Suppose first that $\lambda\subsetneq\mu$. Every cofactor in a column beyond the last changed row vanishes. Indeed, a nonzero minor of an upper triangular matrix requires the deleted row index to be at least the deleted column index. For such a deletion in column $j$, the leading $(j-1)\times(j-1)$ block remains intact and contributes its determinant as a factor. This block contains a zero diagonal entry, so the minor vanishes. All columns supporting $\mathcal T$ lie beyond the last changed row, proving that the discarded tail contributes zero off the diagonal.

Now let $\lambda=\mu$. Entries of $Q$ between blocks of unequal parts vanish. For each maximal block $J$ of equal parts,
\[
 Q_{J,J}=(I_{|J|}+\mathsf S_J)^{-1},\qquad
 Q_{J,J}^{-1}=I_{|J|}+\mathsf S_J,
\]
where $\mathsf S_J$ has $1$ on the first superdiagonal and $0$ elsewhere. In particular, $\det Q=1$, so the tail contribution is
\[
 \operatorname{tr}(Q^{-1}\mathcal T)
 =\sum_{j=1}^n\mathcal T_{jj}
 +\sum_{j=1}^{n-1}\ind_{\lambda_j=\lambda_{j+1}}\mathcal T_{j+1,j}.
\]
The column bound $j\ge n+1-k$ for $\mathcal T$ holds both for the pair $(\lambda,\lambda)$ and for the empty pair with the same $n$ and $k$. For such $j$, we have $j>\ell(\lambda)$, so $a_j=b_j=-j$. If $j<n$, then $a_{j+1}=-(j+1)$ and $\lambda_j=\lambda_{j+1}=0$. Since the discarded sum defining $\mathcal T_{ij}$ depends on the partitions only through $a_i$ and $b_j$, the entries in these two sums and their multipliers agree with those for $\lambda=\mu=\varnothing$. Thus the two pairs have the same tail contribution.

For the empty pair, \eqref{eq:finite} gives zero because $[g_\varnothing]p_k=0$. In the full convolution,
\[
 R^{(k)}_{ii}=\mathsf H_0(k)=0,\qquad
 R^{(k)}_{i+1,i}=\mathsf H_1(k-1)=1.
\]
Thus, with $\mathsf S$ now the $n\times n$ shift matrix,
\[
 \operatorname{tr}\bigl((I+\mathsf S)R^{(k)}\bigr)=n-1.
\]
The discarded tail consequently contributes $n-1$ for the empty pair, and hence for every diagonal pair. Together with the off-diagonal calculation, this proves \eqref{eq:cofactor}.
\end{proof}

Define the formal Laurent series, viewed in $\mathbb Z((z))$,
\begin{equation}\label{eq:uv}
 u_i=z^{-a_i}(1-z)^{-i},\qquad
 v_j=z^{b_j}(1-z)^j,\qquad
 \mathcal K_{\lambda,\mu}^{(n)}(z)=\mathbf v^{\top}\operatorname{adj}(Q)\mathbf u.
\end{equation}
Here $\mathbf u=(u_1,\ldots,u_n)^{\top}$ and $\mathbf v=(v_1,\ldots,v_n)^{\top}$ are column vectors, and ${}^{\top}$ denotes transpose. Since
\[
 [z^k](u_iv_j)=[z^{k+a_i-b_j}](1-z)^{-(i-j)}=R^{(k)}_{ij},
\]
Proposition~\ref{prop:cofactor} and \eqref{eq:variation} give
\begin{equation}\label{eq:Kcoefficient}
 F_k(\lambda,\mu)=[z^k]\mathcal K_{\lambda,\mu}^{(n)}(z)
 -(n-1)\delta_{\lambda\mu}.
\end{equation}
We now evaluate $\mathcal K_{\lambda,\mu}^{(n)}(z)$ exactly as a rational function. For each coefficient $F_k(\lambda,\mu)$, we choose $n$ as in \eqref{eq:nchoice}. The superscript records the matrix size; the off-diagonal expression obtained below is independent of this choice.

\subsection{Geometric evaluation by Schur complements}\label{sec:geometry}
Assume $\lambda\subsetneq\mu$. For $1\le i<n$, define the junction parameter
\[
 \Delta_i=\mu_{i+1}-\lambda_i.
\]
An index $i$ is active if $\mu_i>\lambda_i$. List the active indices as $i_1<\cdots<i_r$, and set
\[
 \mathcal A=\{i_1,\ldots,i_r\},\qquad
 \mathcal I=\{1,\ldots,n\}\setminus\mathcal A,
\]
with both sets ordered increasingly. For ordered index sets $U,V$, write $Q_{U,V}$ for the corresponding submatrix and $\mathbf u_U,\mathbf v_U$ for the corresponding subvectors.

The principal block $Q_{\mathcal I,\mathcal I}$ is unit upper triangular and therefore has determinant $1$. Simultaneously permute the rows and columns to place $\mathcal I$ first; this leaves the determinant unchanged. Define the Schur complement
\[
 \mathcal L=Q_{\mathcal A,\mathcal A}
 -Q_{\mathcal A,\mathcal I}Q_{\mathcal I,\mathcal I}^{-1}Q_{\mathcal I,\mathcal A}.
\]
We index the reduced matrix and vectors by positions $1,\ldots,r$, corresponding to the original indices $i_1,\ldots,i_r$. The symbol $I$, as in $I+\mathsf S$, denotes an identity matrix whose size is determined by the block.

Eliminating inactive indices preserves upper triangularity and leaves the active diagonal entries zero, so $\mathcal L$ is strictly upper triangular. The simultaneous row and column permutations contribute no determinant sign. Applying block elimination to $Q+\varepsilon\mathbf u\mathbf v^{\top}$ and extracting the coefficient of $\varepsilon$ gives
\begin{equation}\label{eq:reduced}
\begin{aligned}
 \mathcal K_{\lambda,\mu}^{(n)}(z)
 &=\mathbf v'{}^{\top}\operatorname{adj}(\mathcal L)\mathbf u',\\
 \mathbf u'&=\mathbf u_{\mathcal A}
 -Q_{\mathcal A,\mathcal I}Q_{\mathcal I,\mathcal I}^{-1}\mathbf u_{\mathcal I},\\
 \mathbf v'{}^{\top}&=\mathbf v_{\mathcal A}^{\top}
 -\mathbf v_{\mathcal I}^{\top}Q_{\mathcal I,\mathcal I}^{-1}Q_{\mathcal I,\mathcal A}.
\end{aligned}
\end{equation}
Indeed, the variation of the inactive block determinant is multiplied by $\det\mathcal L=0$ and therefore vanishes, while its unperturbed determinant is $\det Q_{\mathcal I,\mathcal I}=1$.

\begin{lemma}\label{lem:super}
Let $a=i_p<b=i_{p+1}$ be consecutive active indices. If $b>a+1$, then $\mathcal L_{p,p+1}=0$. If $b=a+1$, then
\[
 \mathcal L_{p,p+1}=\begin{cases}
 1,&\Delta_a=1,\\
 -1,&\Delta_a=0,\\
 0,&\text{otherwise}.
 \end{cases}
\]
\end{lemma}
\begin{proof}
Upper triangularity ensures that only eliminated indices strictly between $a$ and $b$ contribute to this entry. If $b=a+1$, there are no such indices, so
\[
 \mathcal L_{p,p+1}=Q_{a,a+1}
 =\mathsf H_{-1}(1-\Delta_a)=[w^{1-\Delta_a}](1-w),
\]
which gives the asserted values.

Suppose now that $b>a+1$, and put $J=\{a+1,\ldots,b-1\}$. Every index in $J$ is unchanged, and the partition inequalities imply $\mu_b\le\lambda_a$. By \eqref{eq:H}, a nonzero $Q_{a,t}$ with $t\in J$ requires $\lambda_t=\lambda_a$, while a nonzero $Q_{t,b}$ requires $\lambda_t=\mu_b$. Moreover, $Q_{J,J}$ and its inverse have no entries between blocks of unequal part values. Thus, if $\mu_b<\lambda_a$, the Schur-complement correction vanishes. The direct entry $Q_{ab}$ also vanishes, so $\mathcal L_{p,p+1}=0$.

If $\mu_b=\lambda_a$, all intervening part values equal this common value. Hence
\[
 Q_{ab}=(-1)^{b-a},\qquad Q_{J,J}=(I+\mathsf S)^{-1},
\]
where $\mathsf S$ is the first-superdiagonal shift matrix on $J$. Since
\[
 Q_{a,t}=(-1)^{t-a},\qquad Q_{t,b}=(-1)^{b-t}\quad(t\in J),
\]
we have
\[
 (I+\mathsf S)Q_{J,b}=(0,\ldots,0,-1)^{\top}.
\]
Consequently,
\[
 Q_{a,J}Q_{J,J}^{-1}Q_{J,b}=(-1)^{b-a}=Q_{ab}.
\]
The correction therefore cancels the direct entry, completing the proof.
\end{proof}

For $r\ge2$, a strictly upper triangular matrix of order $r$ has only one potentially nonzero minor of order $r-1$, formed from rows $1,\ldots,r-1$ and columns $2,\ldots,r$. Its determinant is the product of the first superdiagonal entries. Thus
\[
 (\operatorname{adj}\mathcal L)_{1r}
 =(-1)^{r-1}\prod_{p=1}^{r-1}\mathcal L_{p,p+1},
\]
and every other entry of $\operatorname{adj}(\mathcal L)$ is zero. The same formula holds for $r=1$ with the empty-product convention, giving $\operatorname{adj}([0])=[1]$.

By Lemma~\ref{lem:super} and \eqref{eq:reduced}, $\mathcal K_{\lambda,\mu}^{(n)}$ can therefore be nonzero only when the active rows are consecutive and $\Delta_i\in\{0,1\}$ at every junction between them. This condition characterizes corner-connected broken ribbons. For two adjacent active rows $i,i+1$, the value $\Delta_i=1$ gives an edge junction, $\Delta_i=0$ gives a corner junction, and $\Delta_i\ge2$ produces a $2\times2$ block. If $\Delta_i<0$, these two rows are separated under edge-or-vertex adjacency. An inactive row between two active rows also separates the diagram.

Suppose now that $\theta$ is a corner-connected broken ribbon with active rows $s,s+1,\ldots,q$, so $r=q-s+1$. There are $r-c$ edge junctions and $c-1$ corner junctions. The unique nonzero entry of $\operatorname{adj}(\mathcal L)$ is therefore
\begin{equation}\label{eq:sign}
 (\operatorname{adj}\mathcal L)_{1r}
 =(-1)^{r-1}\prod_{j=1}^{r-1}\mathcal L_{j,j+1}=(-1)^{r-c}.
\end{equation}
Define the endpoint indicators
\begin{equation}\label{eq:endpoints}
 \tau_-=\begin{cases}
 \ind_{\lambda_{s-1}=\mu_s},&s>1,\\
 0,&s=1,
 \end{cases}
 \qquad\tau_+=\ind_{\lambda_q=\lambda_{q+1}}.
\end{equation}
Since $q$ is the last active row, $\mu_{q+1}=\lambda_{q+1}$ and $\Delta_q=0$ is equivalent to $\tau_+=1$. The effective endpoint vectors in \eqref{eq:reduced} satisfy
\begin{equation}\label{eq:endvectors}
 v'_1=v_s(1-z)^{-\tau_-},\qquad
 u'_r=u_q(1-z)^{-\tau_+}.
\end{equation}
To verify the first identity, note that $v'_1=v_s$ when $s=1$. If $s>1$, all preceding rows are unchanged, so $\mu_s\le\lambda_{s-1}$. Strict inequality makes every prefix entry in column $s$ vanish. If equality holds, only the equal-part prefix block ending at $s-1$ contributes. Multiplying its column by the inverse $I+\mathsf S$ leaves a single nonzero entry, equal to $-1$, at $s-1$. Hence
\[
 v'_1=v_s+v_{s-1}=\frac{v_s}{1-z},\qquad
 \frac{v_{s-1}}{v_s}=\frac{z}{1-z}.
\]
For the second identity, all rows after $q$ are unchanged. If $\lambda_q>\lambda_{q+1}$, every suffix entry in row $q$ vanishes, giving $u'_r=u_q$. If equality holds, multiplying the row on the equal-part suffix block by its inverse $I+\mathsf S$ leaves a single entry $-1$ at $q+1$. Thus
\[
 u'_r=u_q+u_{q+1}=\frac{u_q}{1-z},\qquad
 \frac{u_{q+1}}{u_q}=\frac{z}{1-z}.
\]
Our choice of $n$ ensures that $q+1\le n$. These calculations also cover the case of a single active row.

Equations \eqref{eq:reduced}--\eqref{eq:endvectors} yield
\begin{equation}\label{eq:Kraw}
 \mathcal K_{\lambda,\mu}^{(n)}(z)
 =(-1)^{r-c}z^{\mu_s-\lambda_q+r-1}
 (1-z)^{-(r-1+\tau_-+\tau_+)}.
\end{equation}
The exponents satisfy
\begin{equation}\label{eq:identities}
 d=\mu_s-\lambda_q+(r-c),\qquad B=r-1+\tau_-+\tau_+.
\end{equation}
Indeed, telescoping the row lengths gives
\[
 d=\mu_s-\lambda_q+\sum_{i=s}^{q-1}(\mu_{i+1}-\lambda_i).
\]
Each summand is $1$ at an edge junction and $0$ at a corner junction, so the sum equals $r-c$.

For the second identity, examine the summands in \eqref{eq:boundary}. If $s\le i<q$, then
\[
 \mu_i>\lambda_i\ge\lambda_{i+1},\qquad \lambda_i\le\mu_{i+1},
\]
so the summand equals $1$. The index $s-1$, when present, contributes $\tau_-$, and the index $q$ contributes $\tau_+$. All remaining summands vanish. In particular, $B\ge0$.

Substituting \eqref{eq:identities} into \eqref{eq:Kraw} gives
\begin{equation}\label{eq:Kevaluated}
 \mathcal K_{\lambda,\mu}^{(n)}(z)=(-1)^{r-c}z^{d+c-1}(1-z)^{-B}.
\end{equation}
The right-hand side is independent of $n$ and contains only positive powers of $z$. Thus \eqref{eq:Kcoefficient}, with $n$ chosen as in \eqref{eq:nchoice} for each $k$, proves \eqref{eq:gf}. The formal binomial theorem then gives \eqref{eq:closed}, including $B=0$ under our generalized binomial convention.

For $\lambda=\mu$, the block decomposition by equal parts gives
\[
 \det Q=1,\qquad Q^{-1}=I+\mathsf S_\lambda,
\]
where $(\mathsf S_\lambda)_{i,i+1}=1$ when $\lambda_i=\lambda_{i+1}$, and all other entries are zero. We have $u_iv_i=1$ and, whenever $\lambda_i=\lambda_{i+1}$,
\[
 u_{i+1}v_i=\frac{z}{1-z}.
\]
Since $n>\ell(\lambda)$, all strict drops of $\lambda$ occur among the $n-1$ adjacent pairs. Exactly $n-1-\corn(\lambda)$ of these pairs therefore have equal parts, giving
\begin{equation}\label{eq:Kdiag}
 \mathcal K_{\lambda,\lambda}^{(n)}(z)
 =n+\bigl(n-1-\corn(\lambda)\bigr)\frac{z}{1-z}.
\end{equation}
Every positive-degree coefficient is $n-1-\corn(\lambda)$. Subtracting the correction $n-1$ in \eqref{eq:Kcoefficient} proves \eqref{eq:diagonal}. Together with Lemma~\ref{lem:support} and the preceding off-diagonal calculation, this completes the proof of Theorem~\ref{thm:main}.

\begin{example}\label{ex:two-one-det-proof}
We compute the same entries of Example~\ref{ex:two-one-p-three} as in Example~\ref{ex:two-one-row-proof}, using the Schur complement and \eqref{eq:Kcoefficient}. Take $\lambda=(2,1)$ and $n=11$, which satisfies \eqref{eq:nchoice} for all the cases below.

For $\mu=(3,2)$, the active rows are $1,2$ and $\Delta_1=0$. The reduced matrix and its adjugate are
\[
 \mathcal L=\begin{pmatrix}0&-1\\0&0\end{pmatrix},
 \qquad
 \operatorname{adj}(\mathcal L)=\begin{pmatrix}0&1\\0&0\end{pmatrix}.
\]
Neither endpoint is corrected in \eqref{eq:endvectors}. Thus
\[
 \mathcal K_{\lambda,(3,2)}^{(11)}(z)=v_1u_2
 =z^2(1-z)\cdot z(1-z)^{-2}
 =\frac{z^3}{1-z},
\]
whose coefficient of $z^3$ is $1$.

For $\mu=(2,1,1)$, only row $3$ is active, so $\mathcal L=[0]$ and its adjugate is $[1]$. Both endpoint indicators are $1$, since $\lambda_2=\mu_3=1$ and $\lambda_3=\lambda_4=0$. Using $v_3u_3=z$ in \eqref{eq:endvectors} gives $\mathcal K_{\lambda,(2,1,1)}^{(11)}(z)=z/(1-z)^2$, with coefficient $3$ at $z^3$.

For $\mu=(2,1,1,1)$, the active rows are $3,4$ and $\Delta_3=1$. Now $\mathcal L_{12}=1$, so $(\operatorname{adj}\mathcal L)_{12}=-1$. The two endpoint indicators are again $1$, and $v_3u_4=z^2/(1-z)$. The two endpoint corrections give $\mathcal K_{\lambda,(2,1,1,1)}^{(11)}(z)=-z^2/(1-z)^3$, whose coefficient at $z^3$ is $-3$.

For $\mu=(3,1,1)$, rows $1$ and $3$ are active with an inactive row between them. In this case $Q_{12}=Q_{13}=0$, and the Schur-complement entry $Q_{13}-Q_{12}Q_{22}^{-1}Q_{23}$ vanishes. Thus $\mathcal L$ is the zero matrix of order two, its adjugate vanishes, and the coefficient is $0$.

On the diagonal, the equal-part block of zeros occupies rows $3$ through $11$. There are eight adjacent equal-part pairs, so \eqref{eq:Kdiag} gives $11+8z/(1-z)$. Extracting the coefficient of $z^3$ and subtracting the cutoff correction $n-1=10$ gives $8-10=-2$.
\end{example}

\subsection*{Acknowledgements}
The authors gratefully acknowledge financial support from Da Lat University through research project No.~B2026-TDL-01-KPT.


\begingroup\small
\begin{thebibliography}{99}
\bibitem{AY} A. Amanov and D. Yeliussizov, \emph{Determinantal formulas for dual Grothendieck polynomials}, Proc. Amer. Math. Soc. \textbf{150} (2022), 4113--4128. \href{https://arxiv.org/abs/2003.03907}{arXiv:2003.03907}.
\bibitem{BSZ} J. Bandlow, A. Schilling and M. Zabrocki, \emph{The Murnaghan--Nakayama rule for $k$-Schur functions}, J. Combin. Theory Ser. A \textbf{118} (2011), 1588--1607. \href{https://doi.org/10.1016/j.jcta.2011.01.009}{doi:10.1016/j.jcta.2011.01.009}.
\bibitem{BMS} J. Blasiak, J. Morse and G. H. Seelinger, \emph{$K$-theoretic Catalan functions}, Adv. Math. \textbf{404} (2022), Art.~108421. \href{https://doi.org/10.1016/j.aim.2022.108421}{doi:10.1016/j.aim.2022.108421}; \href{https://arxiv.org/abs/2010.01759}{arXiv:2010.01759}.
\bibitem{Buch} A. S. Buch, \emph{A Littlewood--Richardson rule for the $K$-theory of Grassmannians}, Acta Math. \textbf{189} (2002), 37--78.
\bibitem{GGL} P. Galashin, D. Grinberg and G. Liu, \emph{Refined dual stable Grothendieck polynomials and generalized Bender--Knuth involutions}, Electron. J. Combin. \textbf{23} (2016), no.~3, Paper~P3.14. \href{https://doi.org/10.37236/5737}{doi:10.37236/5737}; \href{https://arxiv.org/abs/1509.03803}{arXiv:1509.03803}.
\bibitem{Iwao} S. Iwao, \emph{Free-fermions and skew stable Grothendieck polynomials}, J. Algebraic Combin. \textbf{56} (2022), 493--526. \href{https://doi.org/10.1007/s10801-022-01121-6}{doi:10.1007/s10801-022-01121-6}; \href{https://arxiv.org/abs/2004.09499v3}{arXiv:2004.09499v3}.
\bibitem{Kim} J. S. Kim, \emph{Jacobi--Trudi formula for refined dual stable Grothendieck polynomials}, J. Combin. Theory Ser. A \textbf{180} (2021), Art.~105415. \href{https://doi.org/10.1016/j.jcta.2021.105415}{doi:10.1016/j.jcta.2021.105415}; \href{https://arxiv.org/abs/2003.00540}{arXiv:2003.00540}.
\bibitem{Kundu} S. Kundu, \emph{A new proof of the Murnaghan--Nakayama-type rule for stable Grothendieck polynomials of Grassmannian type}, Indian J. Pure Appl. Math. (2026). \href{https://doi.org/10.1007/s13226-026-01045-y}{doi:10.1007/s13226-026-01045-y}.
\bibitem{KunduCanonical} S. Kundu, \emph{A note on canonical stable Grothendieck functions}, preprint (2025), \href{https://arxiv.org/abs/2511.00830v1}{arXiv:2511.00830v1}.
\bibitem{LP} T. Lam and P. Pylyavskyy, \emph{Combinatorial Hopf algebras and $K$-homology of Grassmannians}, Int. Math. Res. Not. IMRN (2007), no.~24, Art.~ID rnm125.
\bibitem{LSS} T. Lam, A. Schilling and M. Shimozono, \emph{$K$-theory Schubert calculus of the affine Grassmannian}, Compos. Math. \textbf{146} (2010), 811--852. \href{https://doi.org/10.1112/S0010437X09004539}{doi:10.1112/S0010437X09004539}; \href{https://arxiv.org/abs/0901.1506}{arXiv:0901.1506}.
\bibitem{Macdonald} I. G. Macdonald, \emph{Symmetric Functions and Hall Polynomials}, 2nd ed., Oxford University Press, 1995.
\bibitem{NguyenEtAl} D.-K. Nguyen, D. T. Hiep, T. H. Son and D. L. H. Thuy, \emph{A Murnaghan--Nakayama rule for Grothendieck polynomials of Grassmannian type}, Ann. Comb. \textbf{28} (2024), 155--168.
\bibitem{NguyenKk} D.-K. Nguyen, \emph{A generalization of the Murnaghan--Nakayama rule for $K$-$k$-Schur and $k$-Schur functions}, Int. Math. Res. Not. IMRN \textbf{2024} (2024), no.~6, 4738--4766. \href{https://doi.org/10.1093/imrn/rnad175}{doi:10.1093/imrn/rnad175}; \href{https://arxiv.org/abs/2212.02037v3}{arXiv:2212.02037v3}.
\bibitem{Yeliussizov} D. Yeliussizov, \emph{Duality and deformations of stable Grothendieck polynomials}, J. Algebraic Combin. \textbf{45} (2017), 295--344; \href{https://arxiv.org/abs/1601.01581}{arXiv:1601.01581}.
\bibitem{Yskew} D. Yeliussizov, \emph{Symmetric Grothendieck polynomials, skew Cauchy identities, and dual filtered Young graphs}, J. Combin. Theory Ser. A \textbf{161} (2019), 453--485. \href{https://arxiv.org/abs/1711.09544}{arXiv:1711.09544}.
\end{thebibliography}
\endgroup
\end{document}